\documentclass[conference,a4paper]{IEEEtran}
\IEEEoverridecommandlockouts
\usepackage{cite}
\usepackage{amsmath,amssymb,amsfonts}
\usepackage{algorithmic}

\usepackage{bm}
    
\usepackage{algorithm2e}
 \RestyleAlgo{ruled}    
\usepackage{bbm}
\usepackage{enumerate}

\usepackage{graphicx}
\usepackage{textcomp}
\usepackage{dsfont}
\usepackage{xcolor}
\def\BibTeX{{\rm B\kern-.05em{\sc i\kern-.025em b}\kern-.08em
    T\kern-.1667em\lower.7ex\hbox{E}\kern-.125emX}}

\newcommand{\ignore}[1]{}
    
\newcommand{\X}{{\bf X}}
\newcommand{\rX}{{\bf x}} 

\newcommand{\bU}{{\bar{\bf U}}} 
\newcommand{\rbU}{{\bar{\bf u}}}
\newcommand{\sbU}{\bar{U}}
\newcommand{\rsbU}{\bar{u}}

\renewcommand{\H}{{\bf H}}  
 
\newcommand{\sH}{H}

\newcommand{\sU}{U}

\newcommand{\A}{{\bf A}}  
 
\newcommand{\sA}{A}

\newcommand{\G}{{\bf G}}
\newcommand{\rG}{{\bf g}}
\newcommand{\sG}{G}
\newcommand{\rsG}{g}

\newcommand{\ind}{\mathds{1}}
\newcommand{\bbeta}{{\bm \beta}}

\newtheorem{lemma}{Lemma}
\newcommand{\eop}{{\hfill $\blacksquare$} }
    
\begin{document}

\title{ AoI-Based Opportunistic-Fair mmWave Schedulers\\
\thanks{ The first author's work is partially supported by the Prime Minister's Research Fellowship (PMRF), India. This work is also partially funded by the Technocraft Craft Centre for Applied Artificial Intelligence (TCA2I).}
}

\author{\IEEEauthorblockN{Shiksha Singhal}
\IEEEauthorblockA{\textit{ IEOR, IIT Bombay, India } \\
shiksha.singhal@iitb.ac.in}
\and
\IEEEauthorblockN{Veeraruna Kavitha}
\IEEEauthorblockA{\textit{IEOR, IIT Bombay, India} \\
vkavitha@iitb.ac.in}
\and
\IEEEauthorblockN{Sreenath Ramanath}
\IEEEauthorblockA{\textit{Lekha Wireless, India} 
  \\
sreenath@lekhawireless.com}
}

\maketitle

\begin{abstract}
We consider a system with a  Base Station (BS) and multiple mobile/stationary users. BS uses millimeter waves (mmWaves) for data transmission and hence needs to align beams in the directions of the end-users. The idea is to avail regular user-position estimates, which help in accurate beam alignment towards multiple users, paving way for opportunistic mmWave schedulers. We propose an online algorithm that uses a dual opportunistic and fair scheduler to allocate data as well as position-update channels, in each slot. Towards this, well-known alpha-fair objective functions of utilities of various users,  which further depend upon the age of position-information,  are optimized. We illustrate the advantages of the opportunistic scheduler, by comparing it with the previously proposed mmWave schemes; these schedulers choose one user in each slot and start data transmission only after accurate beam alignment.  We also discuss two ways of introducing fairness in such schemes, both of which perform inferior to the proposed age-based opportunistic scheduler.

\end{abstract}


\vspace{-3mm}
\section{Introduction}
  
With the rapid advancement of technology, the data traffic has increased considerably. Applications like, ultra-high-definition (UHD) 3D video, virtual
and augmented realities, internet-of-things (IoT) etc., demand high data rates. High data rates can be achieved either using high   bandwidth  or transmit power. One cannot increase power arbitrarily  due to health guidelines. To increase bandwidth, the current trend is to use  Millimeter Waves (mmWaves) in the spectral range of 24GHz to 40GHz (\hspace{-.05mm}~\cite{5g}).   

In mmWave communications, the base station (BS) needs to align the beam in the  direction of the end user.  Beam alignment is a challenging task as the location of the users may not be known apriori and further could be varying continuously.  

Opportunistic schedulers (\hspace{-1.4mm}~\cite{tejas,kushner,cellular,debayan}) are widely used in wireless networks to take advantage of a `kind of diversity gain'; the channel conditions are independent across slots and users; in every slot, the BS seeks  channel estimate of each of the users and selects the `best' user for data transmission. The `best' criterion also includes fairness aspects: some far away users have bad channels with high probability and one still needs to allocate a channel to such users to be fair. Generalized $\alpha$-fair opportunistic  schedulers are designed precisely for this purpose: allocate a channel to deprived users when the opportunities are the `best', to an extent depending upon fairness-level defining parameter $\alpha$.

\begin{figure}[!h]
\vspace{-5.5mm}
    \centering
    \includegraphics[width = 0.3\textwidth, height = 0.2\textwidth]{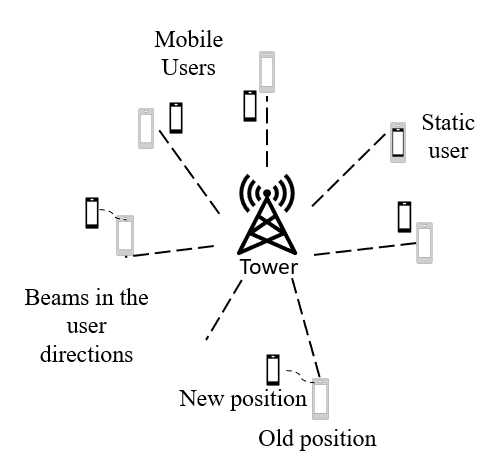}
    \caption{Beam alignment for multiple mobile/static users}
    \label{fig:AoI}
    \vspace{-1mm}
\end{figure}

With mmWave transmissions and its desired accurate beam alignment, it is difficult to obtain the channel estimates from all the users.  Most of the papers (e.g., \cite{ming,qureshi,irmak,mmwave}) that we are aware of select one user in each slot and accurately align the beam towards the selected user prior to data transfer.  Towards designing mmWave-opportunistic scheduler, \textit{we propose to maintain sufficiently accurate user position updates of each user at BS (see Fig. \ref{fig:AoI}).}  With accurate position updates, we assume the beam alignment times to be negligible.
The users can instead transmit the alignment directions, if there are privacy concerns. 

 The wireless users are often mobile,
the BS needs  regular updates of the positions of all the users. Our idea is to design an appropriate dual scheduler that considers all relevant aspects,  age-of-information (AoI), opportunistic and fair schedules;  to be precise, the scheduler has to assign a user for position-update and another/the same user for data transfer in each time-slot in an optimal manner. 
In contrast to the existing literature on AoI that optimize average or peak AoI (e.g., \cite{yuan,kaul,kavitha}), we directly consider the well-known $\alpha$-fair objective function (e.g., \hspace{-1.4mm}~\cite{kushner,cellular}) which is constructed using time-average utilities of all the users, which in turn are influenced by individual AoI trajectories.  

Our main contributions are: i) using theoretical and some heuristic (conjecture that needs proof) arguments, we show that the above mentioned objective function reduces to an average cost Markov decision process (MDP) problem; ii) using the well-known MDP solution techniques (e.g., \cite{Putterman}), we propose an online algorithm for the dual scheduler; iii) the dual scheduler is parametrized by $\alpha$ and  achieves various levels of fairness depending upon the choice of $\alpha$;  iv) we introduce fairness-concepts in two different ways in non-opportunistic schedulers, currently  studied  in literature for mmWaves (e.g., \cite{mmwave}); and finally v) the opportunistic and non-opportunistic schedulers  are compared using exhaustive Monte-Carlo simulations.  


\section{Problem Description and Background}
\label{sec_model}

We consider a network consisting of a Base Station (BS) and $N$ mobile/stationary users labeled $n \in \{1,2,\cdots,N\}$. 
Since millimeter waves (mmWaves) are being used for transmission, one needs to employ narrow beams and a good link is established only when the BS is well-aligned with the users; towards this, we propose to maintain accurate estimates of various user positions at BS; in other words,  we need a mechanism by which regular position updates of (all) the  users is obtained. 


\noindent {\bf User-position Updates:} The time frame is divided into time-slots, and,  exactly one user's position is updated in each slot.  The age of various position-updates  in slot $k$ is represented by   vector $\G_k = (G^1_k, \cdots, G^N_k)$, where $G_k^n = 1$ implies this position is just updated. 
Say user $n$ updates its position in time slot $k$, then the vector $\G_k$ evolves as follows:
\begin{eqnarray}
\sG_{k+1}^n &\hspace{-3mm}=\hspace{-3mm}& 1, \text{ and, } \\
\sG_{k+1}^i &\hspace{-3mm}=\hspace{-3mm}& \min\left\{\sG_{k}^i+1, \ \bar{g}\right\}
\text{ for all } i \neq n, \nonumber
\end{eqnarray}
where $\bar{g}$ is the age upper bound at which the beam alignment deteriorates to a condition such that the utility obtained by that user during data transfer is close to zero.

\noindent
\textbf{Data Scheduler:} The remaining time in every slot after position update is used for data transmission. 
Here again exactly one user is allocated the data channel. The main idea of the paper is to consider an \textit{opportunistic} and \textit{fair scheduler} for this purpose.  Towards opportunistic scheduler, the BS aligns beams in directions of each of the users to obtain their channel estimates, $\H_k = (\sH_k^1,\cdots,\sH_k^n)$; here $\sH^n_k = \sA_k^n \sU_k^n$ in time-slot $k$ for user $n$, where  factor  $\sA_k^n$ depends on the age $\sG^n_k$ of user $n$ and accounts for the misalignment, while,     $\{\sU^n_{k}\}_{k}$ are i.i.d. (independent and identically distributed)  across time-slots for any user $n$; the distribution need not be the same for all users. For example, $\sU^n_{k}$ for any $n$ can depend on Rician or Rayleigh channel conditions of that user.  
If $\sG^n_k$ is high, there is a possibility of higher error in beam alignment as the user might have  moved significantly; then the factor $\sA_k^n$ takes smaller values with higher probability.



\noindent
\textbf{Background on $\mathbf{\alpha}$-Fair Schedulers: }
Fairness is a well-studied concept in wireless networks (e.g., \cite{tejas,kushner,cellular,debayan} and the references therein). Some users are far away from BS while others are nearby. The far away users have inferior channel conditions with high probabilities. Thus, any efficient scheduler (one that maximizes total utility derived because of overall data transfer)  will not be fair to such far away users. 
Fair schedulers are proposed to cater to the demands of such deprived users.

An opportunistic and fair scheduler allocates the channel to the deprived users, in a controlled manner and when the `opportunities'  are the best; this ensures the total utility of the system is the best possible under the given constraints. 
The well-known generalised $\alpha$-fair schedulers exactly achieve this at various levels of fairness indicated  by $\alpha$, by optimizing a certain parameterised concave function of the average utilities obtained by each of the users, as below,
%
\begin{eqnarray}
 \sup_{\bbeta = (\beta^1, \cdots, \beta^N)} \sum_n  \Gamma_\alpha \left ( \rsbU^n(\bbeta) \right ) \text{ with } \rsbU^n(\bbeta):=E[\sH^n \beta^n (\H)]\ ],\hspace{-10mm} \nonumber  \\
 \Gamma_\alpha(\rsbU) := \frac{\rsbU^{1-\alpha} \ind_{ \left \{ \alpha \neq 1 \right \}} }{1-\alpha} + \log(\rsbU) \ind_{\left \{ \alpha =1 \right \}}. 
 \label{Eqn_alpha_fair}
\end{eqnarray}
This criterion was previously considered for i.i.d. channels $\{\H_k\}$, while one can easily extend it to Markov channels (i.e., when $\{\G_k\}$ is Markov). In the following, we consider Markov $\{\G_k\}$,   and, also summarize the solution of the above:
\begin{lemma}
\label{lem_fixed_pt}
{\it Assume $\{ \G_k \}_k$ is  a Markov chain with  at maximum finitely many irreducible classes, each having unique  stationary distribution,  evolving independently of the data scheduler decisions $ \{\bbeta_{k}\}$.   
 Then there exists a unique solution  $(\rsbU_\alpha^1, \cdots, \rsbU_\alpha^N)$ to the following $N$-dimensional fixed point equation:

\vspace{-3mm}
{\small \begin{eqnarray}
\label{Eqn_beta_alpha_FP}
\rsbU_\alpha^n = E^d_{\rX_0}[ \sH^n \beta^n_\alpha (\H) ], \     \beta_\alpha^n (\H) :=\Pi_{i \ne n}  \ind  \left \{ \frac{\sH^n}{ \left (\rsbU_\alpha^n \right )^\alpha }     \ge  \frac{\sH^i}{ \left (\rsbU_\alpha^i \right )^\alpha } \right \}\hspace{-1mm} ,\hspace{-1mm}
\end{eqnarray}}where $E^d_{\rX_0}[\cdot]$ is expectation under stationary distribution(s) and when initialized with $\rX_0$.
Further, $\bbeta_\alpha = (\beta_\alpha^1, \cdots, \beta_\alpha^N)$   optimizes the   $\alpha$-fair criterion \eqref{Eqn_alpha_fair} for any  given $\alpha$.}
\end{lemma}
{\bf Proof:} This is a well known result in literature (e.g., \cite{tejas,cellular} for the case when $\{\H_k\}$ are i.i.d.

  Given the initial condition $\rX_0$ and SMR policy $d$, let $\pi^d_{\rX_0}$ be the stationary distribution (S.D.),  an appropriate convex combination of the S.D.s of various irreducible classes or the unique one. Then one can view $\G$ as i.i.d. variables with this measure, as the given expectation does not depend upon the correlations between various time slots of the Markov chains. And now the proof proceeds as in \cite{tejas,cellular}.
 \eop

\textit{The hypothesis of the above lemma is   satisfied as the age vector $\G_k$ takes finitely many values.}  
An online algorithm that asymptotically represents $\alpha$-fair scheduler \eqref{Eqn_beta_alpha_FP} uses the average utilities  derived by the users till slot $k$, represented by $\bU_k = (\sbU^1_{k},\cdots,\sbU^N_k)$, and is given by (e.g., \cite{kushner,cellular}):

\vspace{-4mm}
{\small
\begin{eqnarray}
\sbU^n_{k+1} &=& \sbU^n_k + \frac{1}{k+1} \left ( \sH^n_{k+1} \beta^n_{k+1} - \sbU^n_k \right ),  \nonumber \\
\beta^n_{k+1} &=& \Pi_{i \ne n}  \ind { \left \{ \frac{\sH^n_{k+1}}{ \left (\sbU^n_k \right )^\alpha }     \ge  \frac{\sH^i_{k+1}}{ \left (\sbU^i_k \right )^\alpha } \right \}  }.  \label{Eqn_beta_scheduler}
\end{eqnarray}}
In \cite{kushner,cellular}, it is proved that the above algorithm converges weakly to asymptotic utilities \eqref{Eqn_beta_alpha_FP} for i.i.d. channels. 
 Observe from \eqref{Eqn_beta_scheduler} that the data scheduler  $\beta_{k+1}$,   depends on the channel estimates $\{\H_{k+1}\}$ as well as the average utilities 
 $\bU_{k}$.
 
\begin{figure}[!h]
     \centering
     \includegraphics[trim = {2cm 8cm 0cm 2cm}, clip, scale = 0.3]{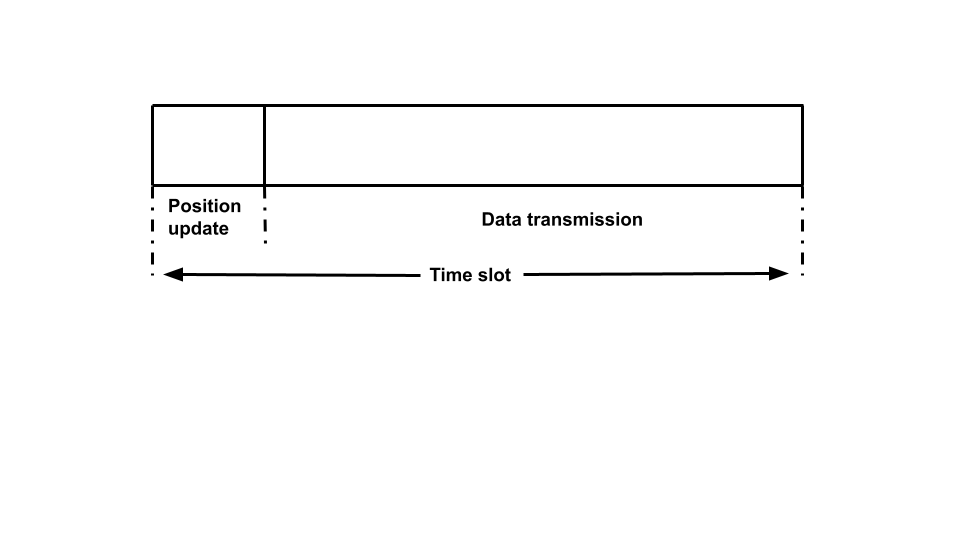}
     \caption{Time slot in proposed mmWave Network}
     \label{fig:time slot division}
 \end{figure}

\noindent 
\textbf{Dual Scheduler: }In view of the above, it is appropriate to consider the age-scheduler that depends upon   $\X_k := \left ( \G_k , \bU_k  \right )$.
Thus to summarise,  we have a dual scheduler in each time slot (see Fig. \ref{fig:time slot division}):  (i) firstly age scheduler
  $d(\X_k)$,  chooses a user  whose position is to be updated, and then, (ii) the data scheduler, $\bbeta(\H_{k+1},\bU_k)$, chooses a user for data transfer.   
  We consider Stationary Markov Randomized (SMR) policies for age-scheduler (e.g., \cite{Putterman}), i.e., $d$ depends only upon $\X_k$. \textit{With slight abuse of notations, we denote any SMR policy by $d$ instead of $d^\infty$.}
  Under such policies, $\X_k$ is a Markov chain which captures the evolution of the entire system.

\section{ 
  Dual Opportunistic Fair Scheduler (DOFS)}
We are interested in solving the following two level optimization problem (if $\bbeta^*$ exists), constructed using average utilities and  $\alpha$-fair objective function defined in \eqref{Eqn_alpha_fair}:
\begin{eqnarray}
\label{Eqn_opt_original}
  \sup_{d \in {\cal D}^1} \sum_{n=1}^N \sbU^n_\infty (\bbeta^*) && \hspace{-10mm}\mbox{ where }
  \bbeta^* \in 
  \arg \max_{ \bbeta} \sum_{n=1}^N  \Gamma_\alpha \left ( \sbU^n_\infty \right )  \mbox{ with }  
  \nonumber 
\\
\sbU^n_\infty &:=&
 \limsup_{T \to \infty} \frac{1}{T}\sum_{k=1}^T  
   \sH_k^n \beta_k^n,
\end{eqnarray}
where the  domain ${\cal D}^1$ includes all the SMR policies.


 
\begin{figure}[!h]
\vspace{-4mm}
    \centering
    \includegraphics[width = 0.4\textwidth, height = 0.2\textwidth]{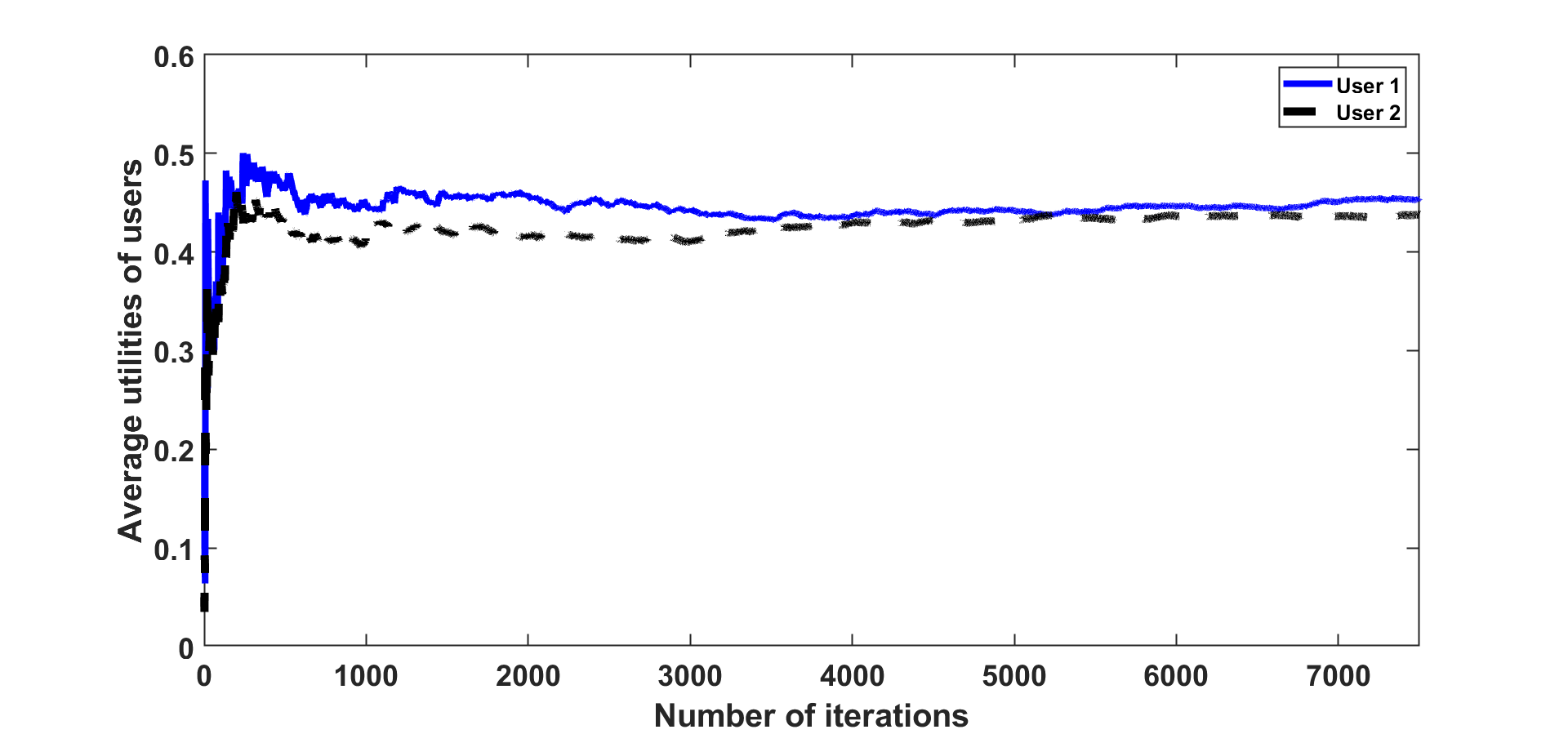}
    \caption{Convergence of average utilities to a constant}
    \label{fig:conjecture}
\end{figure}

 \noindent
{\bf Conjecture:}  For any SMR policy $d$  and any $\bbeta (\H,\bU)$ scheduler, we believe that 
$$
\bU_k (d) \to \rbU^*(d, \bbeta) \mbox{ a.s., }
$$ where $\rbU^*(d, \bbeta)$ is a constant a.s.  This conjecture can be proved either using  Law of Large Numbers for non-homogeneous Markov chains (\hspace{-.05mm}\cite{LLN}) or using stochastic approximation techniques (\hspace{-.05mm}\cite{Benveniste}). We are in the process of constructing this proof, for now we assume the conjecture to be true (see Fig. \ref{fig:conjecture}).

By the above conjecture, 
under any SMR policy $d$ the $\bU_k$, converges to  a constant vector a.s., and hence the SMR policy will start depending only on the first component $\G_k$ of state vector $\X_k$,  and then $\{\H_k\}$ is close to a Markov chain; finally by Lemma \ref{lem_fixed_pt}, the algorithm used in   
\eqref{Eqn_beta_scheduler} converges to the optimal $\alpha$-fair objective function under that given SMR policy. 
Thus the Problem \eqref{Eqn_opt_original}  is equivalent\footnote{All the statements in this paragraph require proofs, we assume these and proceed further to derive online algorithm, which is the aim of this short paper.} to the following, where   the smaller domain $D^2$ includes SMR policies that depend only upon $\G_k$:
\begin{eqnarray}
 \label{Eqn_opt_reduced}
  \sup_{d \in {\cal D}^2} \sum_{n=1}^N \sbU^n_\infty (\bbeta^*) && \hspace{-7
  mm}\mbox{ where }
  \bbeta^* \in 
  \arg \max_{ \bbeta} \sum_{n=1}^N  \Gamma_\alpha \left ( \sbU^n_\infty \right )  \mbox{ with } \nonumber
\\
&&
\hspace{-23mm} \sbU^n_\infty \ := \
 \lim_{T \to \infty} \frac{1}{T}\sum_{k=1}^T   \sH_k^n \beta_k^n \
 \stackrel{a}=  \ E_{\rX_0}^d \left [ \sH^n \beta^n \right ].
\end{eqnarray}
Under any $d \in {\cal D}^2$,  $\{\G_k\}$ itself is a Markov chain,  with finite state space ${\cal S} =  \{ (1, \rsG) \mbox{ or } (\rsG, 1) : 1 < \rsG \le {\bar g}    \}.$ Hence    stationary distribution(s) exist and \underline{hence equality $a$ in} the above; when stationary distribution is unique, the stationary expectation $E^d[.]$ does not depend upon initial condition $\rX_0$. 

By Lemma \ref{lem_fixed_pt}, the inner optimization problem in \eqref{Eqn_opt_reduced} is solved for any $d \in {\cal D}^2$, the solution is given by \eqref{Eqn_beta_alpha_FP} and hence the problem further simplifies, for any initial $\rX_0$, to: 
\begin{eqnarray}
\label{Eqn_P3}
\sup_{d \in {\cal D}^2} \sum_n    \rsbU^n(d)  \hspace{-10mm}&& \mbox{ such that  for each } n, \\
\rsbU^n (d) &=& E^d_{\rX_0} \left[\sH^n  \Pi_{i \ne n}  \ind { \left \{ \frac{\sH^n}{ \left (\rsbU^n (d) \right )^\alpha }     \ge  \frac{\sH^i}{ \left (\rsbU^i (d) \right )^\alpha } \right \}  } \right ]. \nonumber 
\end{eqnarray}
Once again by Lemma \ref{lem_fixed_pt} and Law of Large Numbers for Markov chains (e.g., \cite{LLN}), one can re-write the above as the following average cost MDP problem,

\vspace{-3mm}
{\small \begin{eqnarray}
\label{Eqn_P3_avg}
\sup_{d \in {\cal D}^2} \lim_{T \to \infty} \frac{1}{T}\sum_{k=1}^T E\left [ \sum_n   \sH_{k+1}^n  \Pi_{i \ne n}  \ind { \left \{ \frac{\sH_{k+1}^n}{ \left (\sbU^n_k  \right )^\alpha }     \ge  \frac{\sH_{k+1}^i}{ \left (\sbU^i_k  \right )^\alpha } \right \}  } \right ],
\end{eqnarray}}
where $\{\bU_k\}$
are updated as in \eqref{Eqn_beta_scheduler}.
\begin{algorithm}[h]
\caption{DOFS-initial condition $\rX_0$, $\epsilon$, $\alpha$, $\eta$, $\gamma$}
\label{algo_DOFS}
\begin{enumerate}[(1)]
\item Initialize  $\bU_0$  and $Q(\rX,a)$ for all $\rX,a$ and set $k=0$, $\X_0 = \rX_0$.

\item  {\bf $\epsilon$-greedy age-scheduler:} With probability $\epsilon$,
 \\ random (uniform) age-decision ($d$) is chosen, else \\ choose $d \in \arg \max_{a \in \{1, \cdots, N\}}Q(\X_k,a)$ .

\item Update age of all users, according to the following
\begin{eqnarray*}
&& \sG^n_{k+1} \leftarrow \min \left(\sG^n_{k}+1, \bar{g} \right ) \text{ if } n \ne d \text{ and,}\\
&& \sG^n_{k+1} \leftarrow 1 \text{ for } n = d.
\end{eqnarray*}
\item {\bf Data schedule:} For each $n=1,2,\cdots,N$, set
$$
\beta^n_{k+1}  = \prod_{i \ne n} \ind \left \{ \frac{\sH^n_{k+1}}{(\sbU^n_{k})^\alpha} \geq \frac{\sH^i_{k+1}}{(\sbU^i_{k})^\alpha} \right \}.
$$
\item Update average utilities, for each  $n$,
$$
\sbU^n_{k+1} = \sbU^n_{k}+\frac{1}{k+1}\left(\sH^n_{k+1}\beta^n_{k+1} - \sbU^n_{k} \right).
$$
\item Update \textit{Q-table entry} corresponding to $(\X_k, d)$, \\ using {\small$\X_{k+1} = (\G_{k+1}, \bU_{k+1})$}:

\vspace{-5mm}
{\small 
\begin{eqnarray*}
\hspace{-5mm}
Q(\X_k,d) \ \leftarrow \ Q(\X_k,d) + \eta\bigg (\sum_{n=1}^N \sH^n_{k+1}\beta^n_{k+1} \hspace{-34mm}\\ 
&& \hspace{-15mm} +\gamma \max_a Q(\X_{k+1},a)-Q(\X_k,d)\bigg )
\end{eqnarray*}}
\item Check some convergence criteria for $\bU_k$,
\begin{enumerate}[(1)]
    \item If converged, then stop
    \item Else set $k \leftarrow k+1$ and go to Step 2
\end{enumerate}
\end{enumerate}
\end{algorithm}

\subsection{Online Algorithm-DOFS }

As mentioned before, we have an average cost MDP \eqref{Eqn_P3_avg}. From \cite{AvgtoDiscount},    the optimal policy for average cost MDP can be derived from that of the discounted cost MDP, when discount factor $\gamma$ is close to $1$. By this  observation, we propose an online algorithm which is based on two iterative algorithms: a) the well-known Q-learning algorithm  (e.g., \cite{RL}), and b) the iterative algorithm  that implements \eqref{Eqn_beta_alpha_FP} as in \cite{kushner,cellular}. 

In each time slot, first an age-decision is made using $\epsilon$-greedy scheduler, 
which depends upon the Q-table estimates (as explained in Algorithm \ref{algo_DOFS}). Then an user for  data transfer is chosen using $\bU_k$ and $\H_{k+1}$; the latter are estimated (accurately)  after the user-position updates are obtained. 
Finally all the variables $\bU_{k+1}, \G_{k+1}$ and $Q$-table are  updated.

\section{Non-opportunistic data scheduler}
\label{sec_noops}
In \cite{mmwave}, the authors consider  an optimal user scheduling problem to minimize the beam alignment overhead in mmWave networks, while maintaining the desired QoS (rewards related to data transmission) of each user. In each time slot, the BS selects one user, based on the information related to previous schedules and the beam search algorithm finds the most appropriate beam towards the selected user.
The time spent in aligning the beam, depends upon the  gap between the consecutive slots in which the same user is chosen.
As their purpose is different (maintain QoS, rather than maximize sum utilities as in  Section \ref{sec_model}),  they do not collect channel estimates from all users in any slot.  Thus   their scheduler is non-opportunistic, as  opposed to the one  discussed in Section \ref{sec_model}.  Further more, they do not consider fairness.

We propose two ways of introducing fairness into such schedulers and compare the same with DOFS in Section \ref{sec_examples}. 

\noindent {\bf Dual Non-opportunistic Fair Scheduler (DNOFS):} In each time slot, the scheduler uses the expected conditional channel estimates of the users $\{E[\sH^n_{k+1}|\G_k]\}$  in place of $\{\sH^n_{k+1}\}$,   for data-decision (in \eqref{Eqn_P3_avg}); here we assume $\sH^n_k = L^n_k\sU^n_k$, where $\sU^n_k$ are as before and  the factor $L^n_k$ depends on the age $\sG^n_k$ of user $n$ and   characterises the time lost in aligning the beam (using position update) to the user. The age-scheduler  depends only upon $\G_k$ as in previous section. DNOFS is presented in Algorithm \ref{algo_DNOFS}, which differs from 
Algorithm \ref{algo_DOFS} only in step (4).

\noindent 
{\bf Single Non-opportunistic Fair Scheduler (SNOFS):}
We now consider single decision in each time-slot as in \cite{mmwave}, and include the average utilities so far, $\bU_k$, directly in the objective function to achieve fairness, i.e.,  the age-MDP directly optimizes
the following:

\vspace{-3mm}
{\small \begin{eqnarray}
\label{Eqn_M-DNOFS}
\sup_{d \in {\cal D}^2} \lim_{T \to \infty} \frac{1}{T}\sum_{k=1}^T E\left [ \sum_n  \frac{ L^n_{k+1}(\sG^n_{k+1}) \sU^n_k }{\left (\sbU^n_k  \right )^\alpha} \ind{_{\{n= d(\G_{k+1})\}}}   \right ]
\end{eqnarray}} 
Basically, the data and age scheduler coincide here.

\begin{algorithm}[h]
\caption{DNOFS-initial condition $\rX_0$, $\epsilon$, $\alpha$, $\eta$, $\gamma$}
\label{algo_DNOFS}
Steps (1)-(3) are as in Algorithm \ref{algo_DOFS}.
\begin{enumerate}[(4)]
 
\item {\bf Data schedule:} For each $n=1,2,\cdots,N$, set
$$
\beta^n_{k+1}  = \prod_{i \ne n} \ind \left \{ \frac{E[\sH^n_{k+1}]}{(\sbU^n_k)^\alpha} \geq \frac{E[\sH^i_{k+1}]}{(\sbU^i_k)^\alpha} \right \}
$$

\end{enumerate}
Steps (5)-(7) are as in Algorithm \ref{algo_DOFS}.
 
\end{algorithm}

\section{Numerical Experiments}
\label{sec_examples}

We perform exhaustive Monte Carlo simulations for the proposed DOFS, DNOFS and SNOFS algorithms under different mobility conditions: one where the users are static and other when the users are mobile. 
We consider Rayleigh  channels and hence $\sU^n$ is exponentially distributed with parameter $\lambda^n$. We set $\lambda^1 = 1$, $\lambda^2 = 1.8$ and discount factor, $\gamma = 0.9$. We consider two users, i.e., $N=2$.

\noindent
\textbf{DOFS: } The factor $A_k^n$, which reflects the age dependent channel conditions, is considered to be a binomial random variable as below:
$$
A_k^n = 
\begin{cases}
\bar{a} & \text{with probability } p(\rsG^n_{k}), \\
\underline{a} & \text{else.}
\end{cases}
$$
For high speed users, the probability of good channel $p(\rsG)$ decreases fast with the age, $\rsG$. 
 For the purpose of simulations, we set $\bar{a} = 1.1$ and $\underline{a} = .3$. The values of $p(\rsG)$ can be seen from Fig. \ref{fig:static} and \ref{fig:mobile}. For example, $p(2)$ is 0.9 for static users indicating that the information is not deteriorating fast with time, while for mobile users $p(2) = 0.5$.
 In actuality, the users directly estimate $\sH_k = \sA_k\sU_k$
 when they measure their channels;
 also recall  the estimates are assumed to be accurate.

\noindent
\textbf{DNOFS: }As explained in Section \ref{sec_noops}, the BS  does not collect channel estimates from all the users, it instead utilises the conditional expected channel estimates, $E[\sH^n_k|\G_k]$ of the users as in Algorithm \ref{algo_DNOFS} and recall $\sH^n_k = L^n_k \sU^n_k$.  The factor ${ L^n_k}$ represents the loss owing to the time spent in beam alignment, which depends on the age of the previous position-information. We model it as in the following:
$$
E[L^n_k |\G_k=\rG]  = \frac{\bar{a}p(g^n) + \underline{a}(1-p(g^n))}{\bar{a}},
$$ so that   $E[L^n_k |\G_k] = E[A^n_k |\G_k]$; this is done to ensure fair comparisons. 

\noindent
\textbf{SNOFS: } Here the algorithm considers single scheduling decision, i.e., the same user is selected for   position-update and  data transfer in any time slot. To bring in  fairness, we consider the objective function as in \eqref{Eqn_M-DNOFS}.


\subsection*{Observations}
\begin{figure}[!h]
\vspace{-12mm}
\begin{minipage}{4.cm}
\centering
\includegraphics[trim = {0cm 0cm 0cm 1cm}, clip, scale = 0.18]{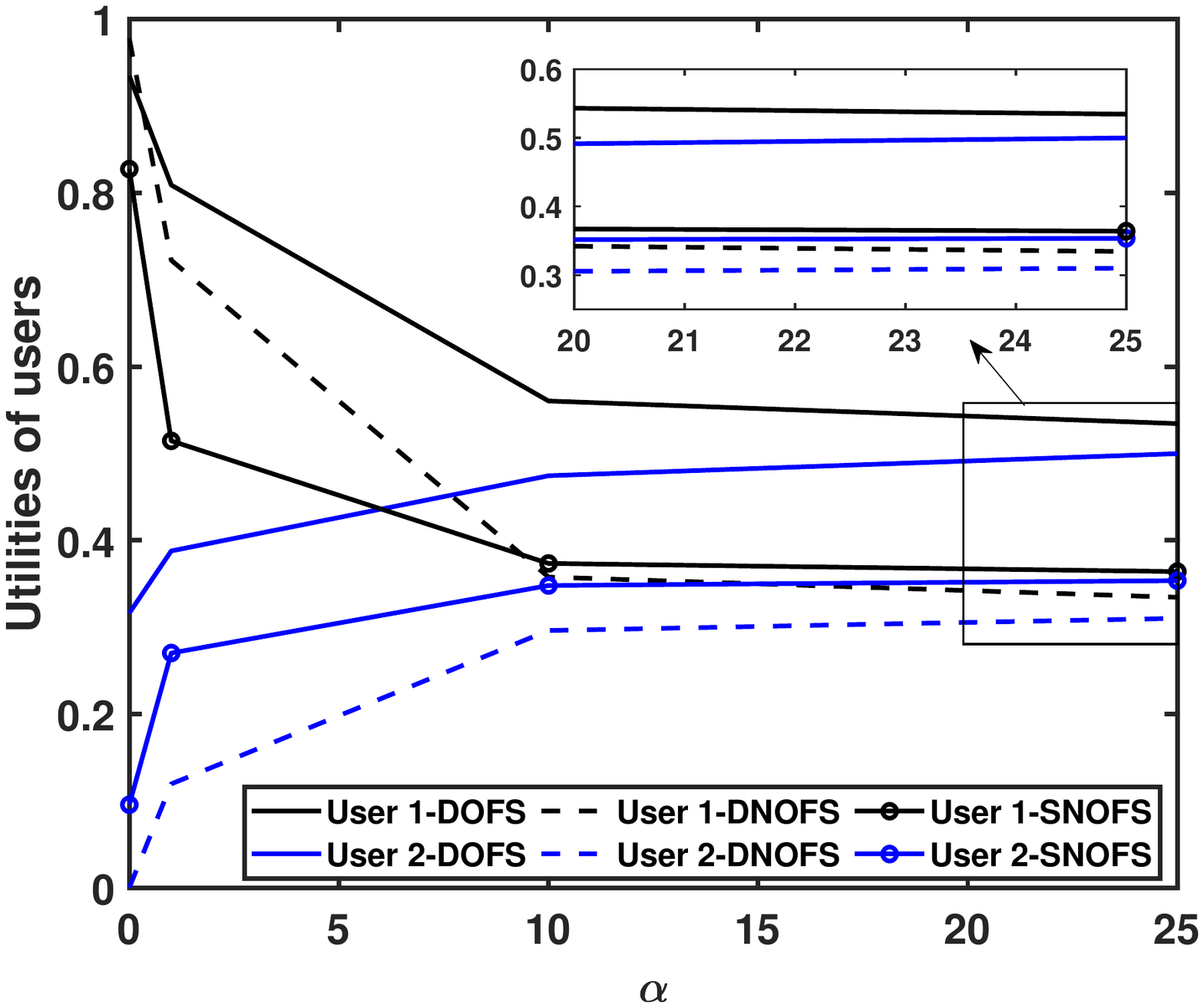}
\end{minipage}\hspace{1.2mm}
\begin{minipage}{4.cm}
\centering
\includegraphics[trim = {0cm 0cm 0cm 1cm}, clip, scale = 0.18]{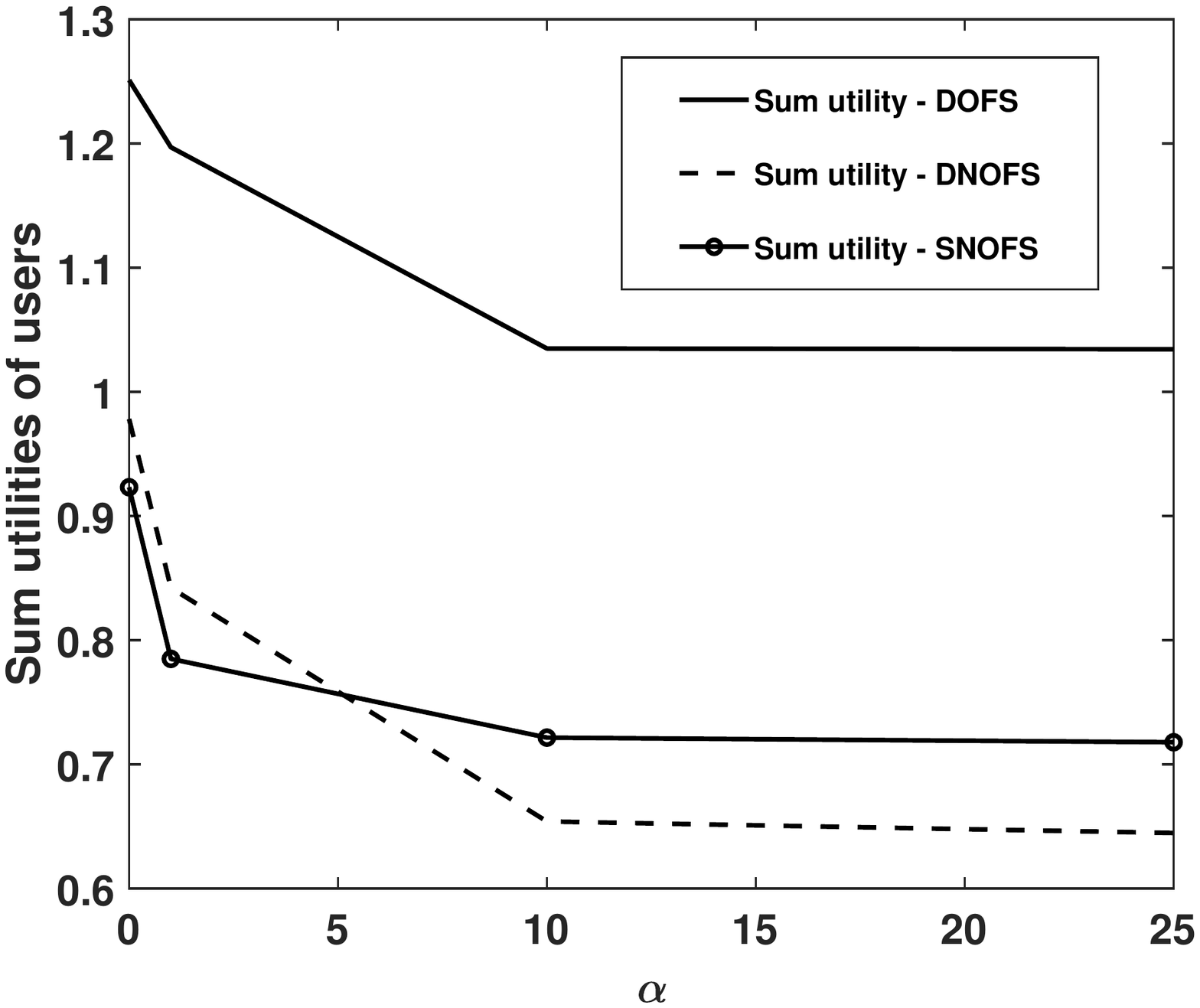}
\end{minipage}
\vspace{-11mm}
\caption{ Utilities (left figure) and sum utilities (right figure) versus fairness parameter $\alpha$ for static users ($p(g) = [1, .9, .5, .1]$).}
\label{fig:static}
\end{figure}

\begin{figure}[!h]
\vspace{-16mm}
\begin{minipage}{4.cm}
\centering
\includegraphics[trim = {0cm 0cm 0cm 1cm}, clip, scale = 0.18]{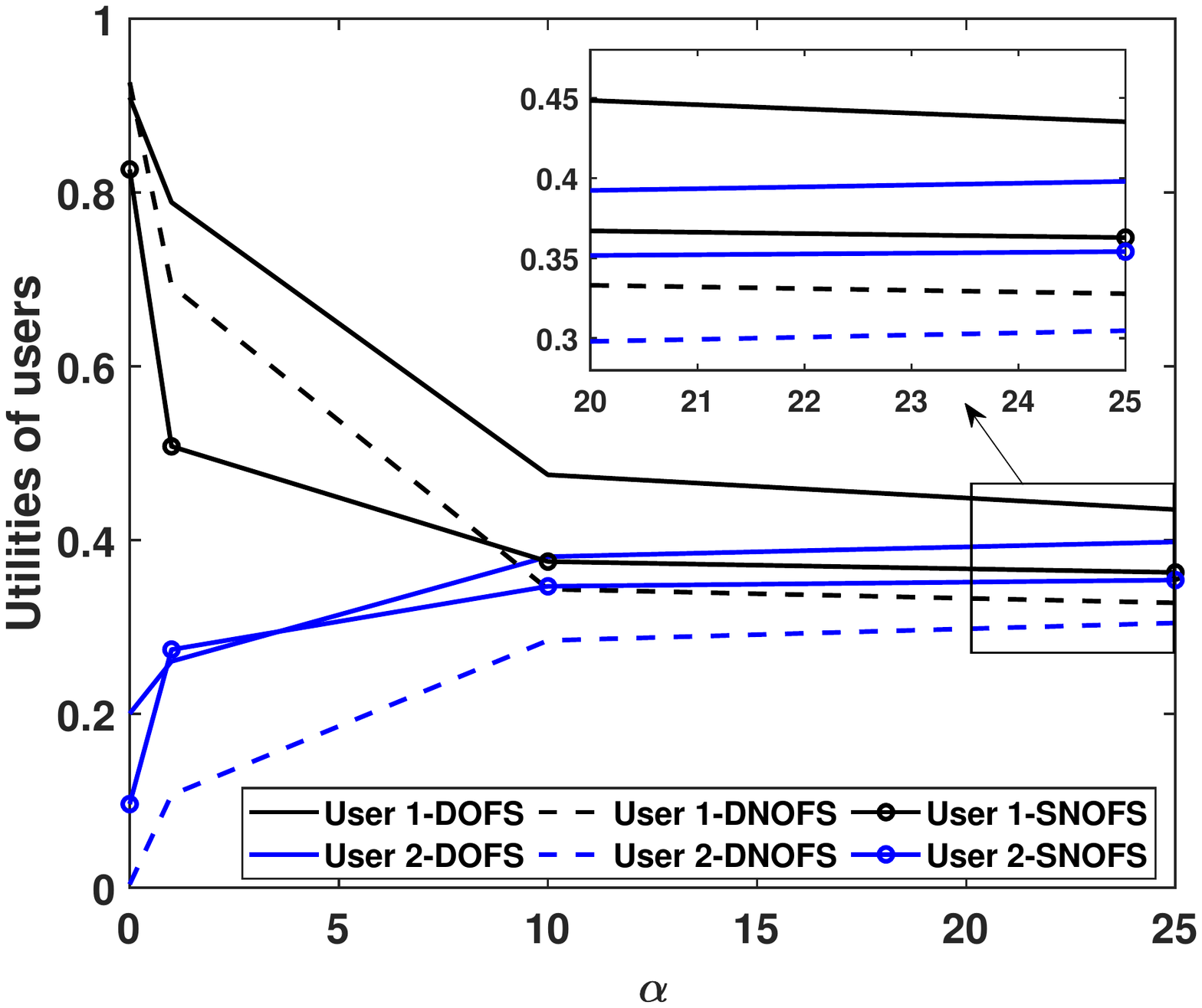}
\end{minipage}
\begin{minipage}{4.cm}
\centering
\includegraphics[trim = {0cm 0cm 0cm 1cm}, clip, scale = 0.18]{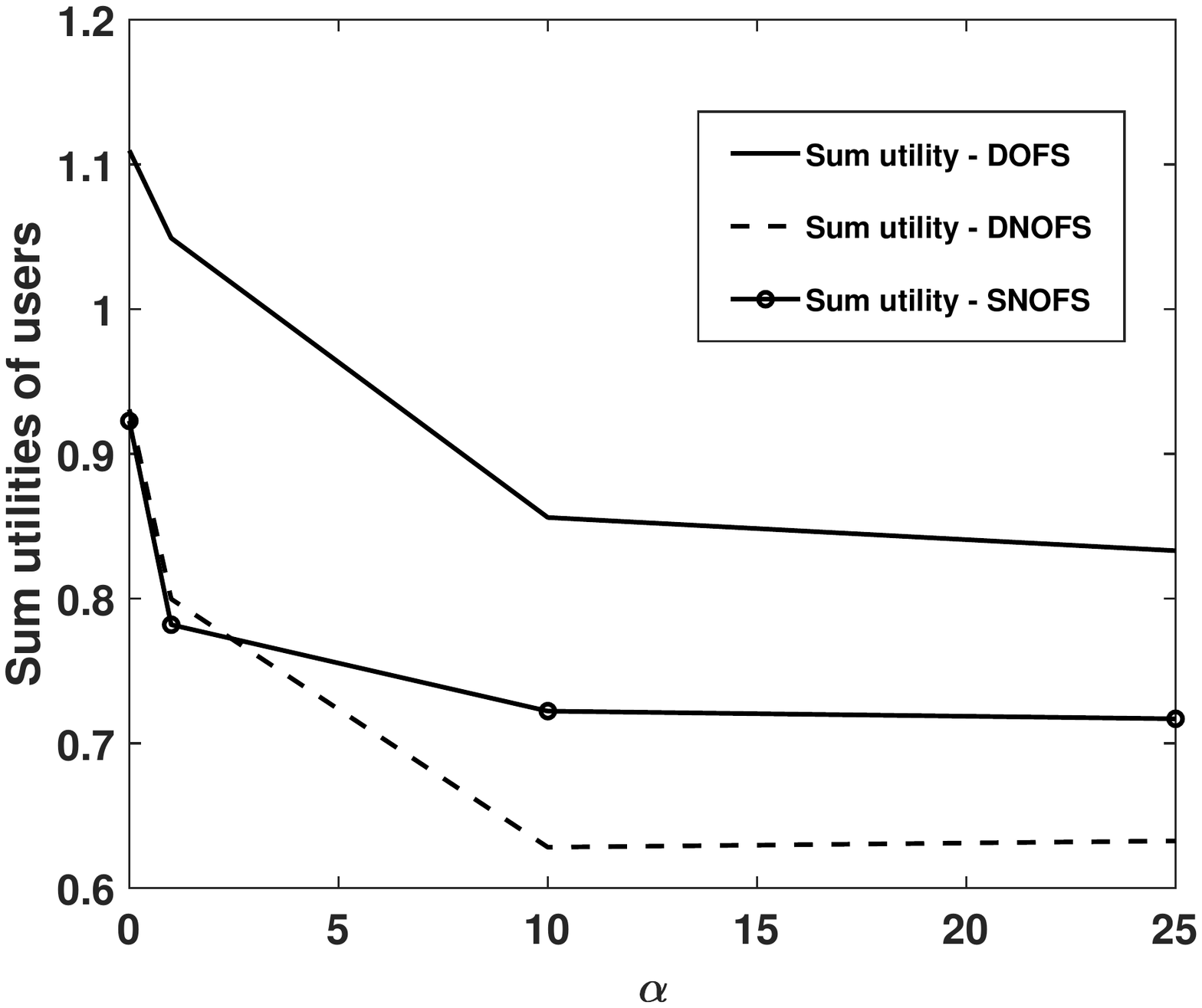}
\end{minipage}
\vspace{-10mm}
\caption{ Utilities (left figure) and sum utilities (right figure) versus fairness parameter $\alpha$ for mobile users ($p(g) = [1, .4, .2, .1]$)}
\label{fig:mobile}
\end{figure}

We analyse the case with $2$ users and observe the following (see Fig. \ref{fig:static} and \ref{fig:mobile}):
(i) The proposed dual scheduler, DOFS outperforms the existing schedulers, DNOFS and SNOFS, for all levels of fairness, $\alpha$; this fact is more clearly evident in the right side figures that have $\sbU^1_\infty+\sbU^2_\infty$ versus $\alpha$; (ii) as the fairness factor $\alpha$ increases, the utilities of the users close in, for all the three schedulers; thus the fair schedulers are effectively ensuring max-min fairness (that equalizes the utilities of all the users) as $\alpha \to \infty.$

\vspace{-3mm}
\section{Conclusions and Future Work}
The millimeter wave communications are the way forward towards supporting the high data-rate applications, as in 5G/6G networks. However, they pose serious design issues; the most complicated issue being the design of narrow and accurate beams directed towards each of the end-users. 

Opportunistic schedulers are known to provide  fair solutions, with a price of fairness that reduces as the number of wireless users increase (\hspace{-.5mm}\cite{mayur}).  Basically, the users with throughout `bad' channels  are allocated the slots when the opportunities  are the `best'  (possible as the channel conditions are independent across users as well as the time slot). 

The opportunistic schedulers require (accurate) estimates of channel conditions of all the users in all the time-slots. With millimeter wave communications that demand accurate beam alignment, one cannot derive channel estimates of all users in every time-slot. We instead propose to maintain a regular and sufficiently accurate estimates of user positions of all the users at the base station. Our inherent assumption is that the beam alignment is possible in negligible time with such accurate position-updates. We achieve this by optimizing the well-known alpha-fair objective functions that further depend upon the age of position-updates of various users.  

We finally propose an online-algorithm that simultaneously implements opportunistic-fair mmWave data-scheduler and an age-scheduler (that updates the user positions). We also propose two methods of incorporating (different levels of) fairness  into 
existing non-opportunistic mmWave-schedulers. We compare the performance of these schedulers with the proposed opportunistic schedulers; the initial simulation results show   that the latter significantly out-performs the former.

This is just the beginning and we have several questions for future investigation. How will  the opportunistic schedulers fare with diverse users, some stationary and some mobile; which type of information (e.g., some estimates of user-speeds) enhances the design of such schedulers; how does the age-scheduler depend upon the statistics of the users; it might be more realistic to consider schedulers over finite-time horizons, how will the age-scheduler depend upon the time average utilities till that slot.

\end{document}